\documentclass[12pt]{article}

\usepackage{latexsym}
\usepackage{amsopn,amscd}
\usepackage{amsfonts}
\usepackage{amssymb}
\usepackage{amsthm}
\usepackage{amsmath}
\usepackage{a4wide}
\usepackage[square,authoryear]{natbib}

\newtheorem{Theorem}{Theorem}[section]

\newtheorem{Lemma}[Theorem]{Lemma}

\def\scs{\sc}
\def\dd{d}

\def\ppartial{\mathcal D}
\def\pppartial{\mathcal D}
\def\ppppartial{[\, \cdot \, , \, \cdot \, ]}
\def\dell{\delta}
\def\ppi{\pi}

\def\fra{\mathfrak}
\def\Sigm{\mathcal S}
\def\jj{j}
\def\nnabla{\widetilde \nabla}
\def\HH{\widetilde H}

\long
\def\MSC#1\EndMSC{\def\arg{#1}\ifx\arg\empty\relax\else
      {\par\narrower\noindent%
      2000 Mathematics Subject Classification: #1\par}\fi}

\long
\def\KEY#1\EndKEY{\def\arg{#1}\ifx\arg\empty\relax\else
    {\par\narrower\noindent%
      Keywords and Phrases: #1\par}\fi}

\def\Nsddata#1#2#3#4#5{
\begin{equation*}
   (#4
     \begin{CD}
      \null @>#2>> \null\\[-3.2ex]
      \null @<<#3< \null
     \end{CD}
    #1, #5)
  \end{equation*}
}

\numberwithin{equation}{section}

\begin{document}

\title{\bf The Lie algebra perturbation Lemma}

\author{
J.~Huebschmann
\\[0.3cm]
 USTL, UFR de Math\'ematiques\\
CNRS-UMR 8524
\\
59655 Villeneuve d'Ascq C\'edex, France\\
Johannes.Huebschmann@math.univ-lille1.fr
 }

\maketitle

\begin{abstract}
{ \noindent We shall establish the following.

\noindent {\bf Theorem.} Let   $(M
     \begin{CD}
      \null @>{\nabla}>> \null\\[-3.2ex]
      \null @<<{\pi}< \null
     \end{CD}
    \fra g, h)$
be a {\em contraction\/} of chain complexes and suppose that {\sl
$\fra g$ is endowed with a bracket $\ppppartial$ turning it into
differential graded Lie algebra. Then the given contraction and
the bracket $\ppppartial$ determine an sh-Lie algebra structure on
$M$, that is, a coalgebra perturbation $\ppartial$ of the
coalgebra differential $\dd^0$ on (the cofree coaugmented
differential graded cocommutative coalgebra) $\Sigm^{\mathrm
c}[sM]$ (on the suspension $sM$ of $M$), the coalgebra
differential $\dd^0$ being induced by the differential on $M$, a
Lie algebra twisting cochain $\tau \colon
\Sigm_{\ppartial}^{\mathrm c}[sM] \longrightarrow \fra g$ and,
furthermore, a contraction
\begin{equation*}
   \left(\Sigm_{\ppartial}^{\mathrm c}[sM]
     \begin{CD}
      \null @>{\overline \tau}>> \null\\[-3.2ex]
      \null @<<{\Pi}< \null
     \end{CD}
    \mathcal C[\fra g], H\right)
  \end{equation*}
of chain complexes which are natural in terms of the data.}

\noindent
Here  $\mathcal C[\fra g]$ refers to the classifying
coalgebra of $\fra g$. }
\end{abstract}

\section{Introduction}
The purpose of this paper is to establish the {\em Perturbation
Lemma\/} for differential graded Lie algebras. The main technique
is H(omological)P(erturbation)T(he\-o\-ry).

A special case of the Lie algebra perturbation Lemma has been
explored in \cite{huebstas} (Theorem 2.7), and the details of the
proof of that Theorem have been promised to be given elsewhere;
the present paper includes these details. The main application of
the quoted Theorem in \cite{huebstas} was the construction of
solutions of the {\em master equation\/} under suitable general
circumstances, and the present paper in particular yields
solutions of the master equation under circumstances even more
general than those in \cite{huebstas}. In the present paper, we
will not elaborate on the master equation, though; suffice it to
mention that the master equation amounts to the defining equation
of a {\em Lie algebra twisting cochain\/} which, in turn, will be
reproduced as \eqref{master} below. Detailed comments related with
the master equation may be found in \cite{huebstas}.

The ordinary perturbation Lemma for chain complexes (reproduced
below as Lemma \ref{ordinary}) has become a standard tool to
handle higher homotopies in a constructive manner. This Lemma is
somehow lurking behind the formulas (1) in Ch. II \S1 of
\cite{shih}, seems to have first been made explicit by {M.
Barratt} (unpublished) and, to our knowledge, appeared first in
print in \cite{rbrown}. Thereafter it has been exploited at
various places in the literature, cf. among others
\cite{gugenhtw}--\cite{huebstas}. The basic reason why HPT works
is the old observation that an exact sequence of chain complexes
which splits as an exact sequence of graded modules and which has
a contractible quotient necessarily splits in the category of
chain complexes \cite{dold} (2.18).

Some more historical comments about HPT may be found e.~g. in
\cite{jimmurra} and in
Section 1 (p. 248) and Section 2 (p. 261) of \cite{huebkade},
which has one of the strongest results in relation to {\em
compatibility with other such as algebra or coalgebra
structure\/}, since it was perhaps first recognized in
\cite{perturba}. Suitable HPT constructions that are compatible
with other algebraic structure enabled us to carry out complete
numerical calculations in group cohomology \cite{cohomolo},
\cite{modpcoho}, \cite{intecoho} which cannot be done by other
methods.

In view of the result of Kontsevich's that the Hochschild complex
of the algebra of smooth functions on a smooth manifold, endowed
with the Gerstenhaber bracket, is formal as a differential graded
Lie algebra \cite{kontssev}, sh-Lie algebras have become a
fashionable topic. The attempt to treat the corresponding higher
homotopies by means of  a suitable version of HPT, relative to the
{\em requisite additional algebraic structure\/}, that is, to make
the perturbations {\em compatible with Lie brackets or more
generally with sh-Lie structures\/}, led to the paper
\cite{huebstas}, but technical complications arise since
 the {\em tensor trick\/}, which was successfully exploited in
\cite{gulstatw} and \cite{huebkade}, breaks down for cocommutative
coalgebras; indeed, the notion of {\em homotopy of morphisms of
cocommutative coalgebras is a subtle concept\/} \cite{schlstas},
and only a special case was handled in \cite{huebstas}, with some
of the technical details merely sketched. The present paper
provides {\em all\/} the necessary details and handles the case of
a {\em general contraction\/} whereas in \cite{huebstas} only the
case of a contraction of a differential graded Lie algebra onto
its {\em homology\/} was treated.

In a subsequent paper \cite{pertlie2} we have extended the
perturbation Lemma to the more general situation of sh-Lie
algebras.

I am much indebted to Jim Stasheff for having prodded me on
various occasions to pin down the perturbation Lemma for Lie
algebras as well as for a number of comments on a draft of the
paper, and to J. Grabowski and P. Urbanski for discussions about
the symmetric coalgebra.

\section{The Lie algebra perturbation Lemma}
\label{outline}

To spell  out the Lie algebra perturbation Lemma, and to
illuminate the unexplained terms in the introduction, we need some
preparation.

The ground ring is a commutative ring with $1$ and will be denoted
by $R$. Perhaps some condition has to be imposed upon $R$ but $R$
is not necessarily a field. In Section \ref{pre} below we will
make this precise. We will take {\em chain complex\/} to mean {\em
differential graded\/} $R$-{\em module\/}. A chain complex will
not necessarily be concentrated in non-negative or non-positive
degrees. The differential of a chain complex will always be
supposed to be of degree $-1$. For a filtered chain complex $X$, a
{\em perturbation\/} of the differential $d$ of $X$ is a
(homogeneous) morphism $\partial$ of the same degree as $d$ such
that $\partial$ lowers the filtration and $(d +
\partial)^2 = 0$ or, equivalently,
\begin{equation}
[d,\partial] + \partial \partial = 0.
\end{equation}
Thus, when $\partial$ is a perturbation on $X$, the sum $d +
\partial$, referred to as the {\em perturbed differential\/},
endows $X$ with a new differential. When $X$ has a graded
coalgebra structure such that $(X,d)$ is a differential graded
coalgebra, and when the {\em perturbed differential\/} $d +
\partial$ is compatible with the graded coalgebra structure, we
refer to $\partial$ as a {\em coalgebra perturbation\/}; the
notion of {\em algebra perturbation\/} is defined similarly. Given
a differential graded coalgebra $C$ and a coalgebra perturbation
$\partial$ of the differential $d$ on $C$, we will occasionally
denote the new or {\em perturbed\/} differential graded coalgebra
by $C_{\partial}$.

A {\em contraction\/}
\begin{equation}
   (N
     \begin{CD}
      \null @>{\nabla}>> \null\\[-3.2ex]
      \null @<<{\pi}< \null
     \end{CD}
    M, h) \label{co}
  \end{equation}
of chain complexes \cite{eilmactw} consists of

-- chain complexes $N$ and $M$,
\newline
\indent -- chain maps $\pi\colon N \to M$ and $\nabla \colon M \to
N$,
\newline
\indent --  a morphism $h\colon N \to N$ of the underlying graded
modules of degree 1;
\newline
\noindent these data are required to satisfy
\begin{align}
 \pi \nabla &= \mathrm{Id},
\label{co0}
\\
Dh &= \mathrm{Id} -\nabla \pi, \label{co1}
\\
\pi h &= 0, \quad h \nabla = 0,\quad hh = 0. \label{side}
\end{align}
The requirements \eqref{side} are referred to as {\em annihilation
properties\/} or {\em side conditions\/}.

Let $\fra g$ be (at first) a chain complex, the differential being
written as $d \colon \fra g \to \fra g$, and let
\begin{equation}
   (M
     \begin{CD}
      \null @>{\nabla}>> \null\\[-3.2ex]
      \null @<<{\pi}< \null
     \end{CD}
    \fra g, h)
\label{cont1}
  \end{equation}
be a {\em contraction\/} of chain complexes; later we will take
$\fra g$ to be a  differential graded Lie algebra. In the special
case where the differential on $M$ is zero, $M$ plainly amounts to
the homology $\mathrm H(\fra g)$ of $\fra g$; in this case, with
the notation $\mathcal H = \nabla \mathrm H(\fra g)$, the
resulting decomposition
\[
\fra g = d \fra g \oplus \mathrm{ker} (h)  =d \fra g \oplus
\mathcal H \oplus h \fra g
\]
may be viewed as a generalization of the familiar {\em Hodge\/}
decomposition.

Let $C$ be a {\em coaugmented\/} differential graded coalgebra
with coaugmentation map $\eta \colon R \to C$ and {\em
coaugmentation\/} coideal $JC = \mathrm{coker}(\eta)$, the
diagonal map being written as $\Delta \colon C \to C \otimes C$ as
usual. Recall that the counit $\varepsilon \colon C \to R$ and the
coaugmentation map determine a direct sum decomposition $C = R
\oplus JC$. The {\em coaugmentation\/} filtration $\{\mathrm
F_nC\}_{n \geq 0}$ is as usual given by
\[\mathrm F_nC = \mathrm{ker}(C \longrightarrow (JC)^{\otimes
(n+1)})\  (n \geq 0)
\]
where the unlabelled arrow is induced by some iterate of the
diagonal $\Delta$ of $C$. This filtration is well known to turn
$C$ into a {\em filtered\/} coaugmented differential graded
coalgebra; thus, in particular, $\mathrm F_0C = R$. We recall that
$C$ is said to be {\em cocomplete\/} when $C=\cup \mathrm F_nC$.

Write $s$ for the {\em suspension\/} operator as usual and
accordingly $s^{-1}$ for the {\em desuspension\/} operator. Thus,
given the chain complex $X$, $(sX)_j = X_{j-1}$, etc., and the
differential ${d\colon sX \to sX}$ on the suspended object $sX$ is
defined in the standard manner so that $ds+sd=0$. Let
$\Sigm^{\mathrm c} = \Sigm^{\mathrm c}[sM]$, the {\em cofree\/}
coaugmented differential graded {\em cocommutative\/} coalgebra
or, equivalently, differential graded {\em symmetric\/} coalgebra,
on the suspension $sM$ of $M$. This kind of coalgebra is well
known to be cocomplete; the existence of its diagonal map may
require some mild assumptions which we will comment upon in
Section \ref{pre} below---requiring $M$ to be projective as an
$R$-module or requiring $R$ to contain the field of rational
numbers as a subring will certainly suffice. Further, let
$\dd^0\colon \Sigm^{\mathrm c} \longrightarrow \Sigm^{\mathrm c}$
denote the coalgebra differential on $\Sigm^{\mathrm c} =
\Sigm^{\mathrm c}[sM]$ induced by the differential on $M$. For $b
\geq 0$, we will henceforth denote the homogeneous degree $b$
component of $\Sigm^{\mathrm c}[sM]$ by $\Sigm_b^{\mathrm c}$;
thus, as a chain complex, $\mathrm F_b\Sigm^{\mathrm c} = R \oplus
\Sigm_1^{\mathrm c} \oplus \dots \oplus \Sigm_b^{\mathrm c}$.
Likewise, as a chain complex, $\Sigm^{\mathrm c} =
\oplus_{j=0}^{\infty} \Sigm_j^{\mathrm c}$. We denote by
\[
\tau_{M}\colon \Sigm^{\mathrm c} \longrightarrow M
\]
the composite of the canonical projection $\mathrm{proj}\colon
\Sigm^{\mathrm c} \to sM$ from $\Sigm^{\mathrm c} = \Sigm^{\mathrm
c}[sM]$ to its homogeneous degree 1 constituent $sM$ with the
desuspension map $s^{-1}$ from $sM$ to $M$.

Given two chain complexes $X$ and $Y$, recall that
$\mathrm{Hom}(X,Y)$ inherits the structure of a chain complex by
the operator $D$ defined by
\begin{equation}
D \phi = d \phi -(-1)^{|\phi|} \phi d
\end{equation}
where $\phi$ is a homogeneous homomorphism from $X$ to $Y$ and
where $|\phi|$ refers to the degree of $\phi$.

Consider the cofree coaugmented differential graded cocommutative
coalgebra $\Sigm^{\mathrm c}[s\fra g]$ on the suspension $s\fra g$
of $\fra g$ (the existence of which we suppose) and, as before,
let
\[
\tau_{\fra g}\colon \Sigm^{\mathrm c}[s\fra g]\longrightarrow \fra
g
\]
be the composite of the canonical projection to $\Sigm_1^{\mathrm
c}[s\fra g] =s\fra g$ with the desuspension map. Suppose that
$\fra g$ is endowed with a graded skew-symmetric bracket $[\,
\cdot \, , \, \cdot \, ]$ that is compatible with the differential
but not necessarily a graded Lie bracket, i.~e. does not
necessarily satisfy the graded Jacobi identity. Let $C$ be a
coaugmented differential graded cocommutative coalgebra. Given
homogeneous morphisms $a,b \colon C \to \fra g$, with a slight
abuse of the bracket notation $[\, \cdot \, , \, \cdot \, ]$, the
{\em cup bracket\/} $[a, b]$ is given by the composite
\begin{equation}
\begin{CD} C @>{\Delta}>> C\otimes C @>{a\otimes b}>>\fra g
\otimes\fra g @> {[\cdot,\cdot]}>> \fra g.
\end{CD}
\label{cupbracket}
\end{equation}
The cup bracket $[\, \cdot \, , \, \cdot \, ]$ is well known to be
a graded skew-symmetric bracket on $\mathrm{Hom}(C,\fra g)$ which
is compatible with the differential on $\mathrm{Hom}(C,\fra g)$.
Define the coderivation
\[
\partial\colon\Sigm^{\mathrm
c}[s\fra g] \longrightarrow \Sigm^{\mathrm c}[s\fra g]
\]
on $\Sigm^{\mathrm c}[s\fra g]$  by the requirement
\begin{equation}
\tau_{\fra g} \partial = \tfrac 12 [\tau_{\fra g},\tau_{\fra g}]
\colon \Sigm_2^{\mathrm c}[s\fra g] \to \fra g. \label{proc333}
\end{equation}
Then $D\partial\  (=d\partial + \partial d) = 0$ since the bracket
on $\fra g$ is supposed to be compatible with the differential
$d$. Moreover, {\sl the bracket on $\fra g$ satisfies the graded
Jacobi identity if and only if $\partial\partial = 0$, that is, if
and only if $\partial$ is a coalgebra perturbation of the
differential $d$ on\/} $\Sigm^{\mathrm c}[s\fra g]$, cf. e.~g.
\cite{huebstas}. The Lie algebra perturbation Lemma below will
generalize this observation.

We now suppose that the graded bracket $[\, \cdot \, , \, \cdot \,
]$ on $\fra g$ turns $\fra g$  into a differential graded Lie
algebra and continue to denote the resulting coalgebra
perturbation by $\partial$, so that $\Sigm_{\partial}^{\mathrm
c}[s\fra g]$ is a coaugmented differential graded cocommutative
coalgebra; in fact, $\Sigm_{\partial}^{\mathrm c}[s\fra g]$ is
then precisely the ordinary {\scs
C(artan-)C(hevalley-)E(ilenberg)\/} or {\em classifying\/}
coalgebra for $\fra g$ and, following \cite{quilltwo} (p.~291), we
denote it by $\mathcal C[\fra g]$ (but the construction given
above is different from that in \cite{quilltwo} which, in turn, is
carried out only over a field of characteristic zero).
Furthermore, given a coaugmented differential graded cocommutative
coalgebra $C$, the cup bracket turns $\mathrm{Hom}(C,\fra g)$ into
a differential graded Lie algebra. In particular,
$\mathrm{Hom}(\Sigm^{\mathrm c},\fra g)$ and $\mathrm{Hom}(\mathrm
F_n\Sigm^{\mathrm c},\fra g)$ ($n\geq 0$) acquire differential
graded Lie algebra structures.

Given a  coaugmented differential graded cocommutative coalgebra
$C$ and a differential graded Lie algebra $\fra h$, a {\em Lie
algebra twisting cochain\/} $t \colon C \to \fra h$ is a
homogeneous morphism of degree $-1$ whose composite with the
coaugmentation  map is zero and which satisfies
\begin{equation}
Dt = \tfrac 12 [t,t] , \label{master}
\end{equation} cf.
\cite{moorefiv}, \cite{quilltwo}. In particular, relative to the
graded Lie bracket $\ppppartial$ on $\fra g$, $\tau_{\fra g}\colon
\mathcal C[\fra g] \to \fra g$ is a Lie algebra twisting cochain,
the {\scs C(artan-)C(hevalley-)E(ilenberg)\/} or {\em universal\/}
Lie algebra twisting cochain for $\fra g$. It is, perhaps, worth
noting that, when $\fra g$ is viewed as an abelian differential
graded Lie algebra relative to the zero bracket, $\Sigm^{\mathrm
c}[s\fra g]$ is the corresponding {\scs CCE\/} or {\em
classifying\/} coalgebra and $\tau_{\fra g}\colon \Sigm^{\mathrm
c}[s\fra g] \to \fra g$ is still the universal Lie algebra
twisting cochain. Likewise, when $M$ is viewed as an {\em
abelian\/} differential graded Lie algebra, $\Sigm^{\mathrm
c}=\Sigm^{\mathrm c}[sM]$ may be viewed as the {\scs CCE\/} or
{\em classifying\/} coalgebra $\mathcal C [M]$ for $M$, and
$\tau_{M} \colon \Sigm^{\mathrm c} \to M$ is then the universal
differential graded Lie algebra twisting cochain for $M$.

At the risk of making a mountain out of a molehill, we note that,
in \eqref{proc333} and \eqref{master} above, the factor $\tfrac
12$ is a mere matter of convenience. The correct way of phrasing
graded Lie algebras when the prime 2 is not invertible in the
ground ring is in terms of an additional operation, the {\em
squaring\/} operation $\mathrm {Sq}\colon \fra g_{\mathrm{odd}}
\to \fra g_{\mathrm{even}}$ and, by means of this operation, the
factor $\tfrac 12$ can be avoided. Indeed, in terms of this
operation, the equation \eqref{master} takes the form
\[
Dt=\mathrm{Sq}(t).
\]
For intelligibility, we will follow the standard convention, avoid
spelling out the squaring operation explicitly, and keep the
factor $\tfrac 12$. A detailed description of the requisite
modofications when the prime 2 is not invertible in the ground is
given in \cite{pertlie2}.

Given a chain complex $\fra h$, an sh-{\em Lie algebra
structure\/} or $L_{\infty}$-{\it structure\/} on $\fra h$ is a
{\em coalgebra perturbation\/} $\partial$ of the differential $d$
on the coaugmented differential graded cocommutative coalgebra
$\Sigm^{\mathrm c}[s \fra h]$ on $s \fra h$, cf. \cite{huebstas}
(Def. 2.6). Given two sh-Lie algebras $(\fra h_1, \partial_1)$ and
$(\fra h_2,\partial_2)$, an {\em sh-morphism\/} or {\em sh-Lie
map\/} from $(\fra h_1, \partial_1)$ to $(\fra h_2, \partial_2)$
is a morphism ${ \Sigm^{\mathrm c}_{\partial_1}[s\fra h_1] \to
\Sigm^{\mathrm c}_{\partial_2}[s\fra h_2] }$ of coaugmented
differential graded coalgebras, cf. \cite{huebstas}.

\begin{Theorem}[Lie algebra perturbation Lemma]
\label{lem11} Suppose that $\fra g$ carries a differential graded
Lie algebra structure. Then the contraction {\rm {\eqref{cont1}}}
and the graded Lie algebra structure on $\fra g$ determine an
sh-Lie algebra structure on $M$, that is, a coalgebra perturbation
$\ppartial$ of the coalgebra differential $\dd^0$ on
$\Sigm^{\mathrm c}[sM]$, a Lie algebra twisting cochain
\begin{equation}
\tau \colon \Sigm_{\ppartial}^{\mathrm c}[sM] \longrightarrow \fra
g \label{tc1}
\end{equation}
and, furthermore, a contraction
\begin{equation}
   \left(\Sigm_{\ppartial}^{\mathrm c}[sM]
     \begin{CD}
      \null @>{\overline \tau}>> \null\\[-3.2ex]
      \null @<<{\Pi}< \null
     \end{CD}
    \mathcal C[\fra g], H\right)
\label{cont5}
  \end{equation}
of chain complexes which are natural in terms of the data so that
\begin{align}
\pi \tau&=\tau_M\colon \Sigm^{\mathrm c}[sM] \longrightarrow M,
\label{twist33}
\\
h \tau &= 0, \label{twist44}
\end{align}
and so that, since by construction, the {\em injection\/}
$\overline \tau \colon \Sigm_{\ppartial}^{\mathrm c}[sM]\to
\mathcal C[\fra g]$ of the contraction is the adjoint $\overline
\tau$ of $\tau$, this injection is then a {\em morphism\/} of
coaugmented differential graded coalgebras.
\end{Theorem}

In the statement of this theorem, the perturbation $\ppartial$
then encapsulates the asserted sh-Lie structure on $M$, and the
adjoint $\overline \tau$ of \eqref{tc1} is plainly an
sh-equivalence in the sense that it induces an isomorphism on
homology, including the brackets of all order that are induced on
homology.

The proof of Theorem \ref{lem11} to be given below includes, in
particular, a proof of Theorem 2.7 in \cite{huebstas}; in fact,
the statement of that theorem is the special case of Theorem
\ref{lem11} where the differential on $M$ is zero, and the details
of the proof of that theorem had been promised to be given
elsewhere.

Theorem \ref{lem11} asserts not only the existence of the Lie
algebra twisting cochain \eqref{tc1} and contraction \eqref{cont5}
but also includes explicit natural constructions for them, under
the additional assumption that the prime 2 be invertible in the
ground ring. The necessary modifications for the general case
where the prime 2 is not necessarily invertible in the ground ring
are explained in \cite{pertlie2}. The explicit constructions for
the coalgebra perturbation $\ppartial$ and Lie algebra twisting
cochain \eqref{tc1} will be spelled out in Complement I below, and
explicit constructions of the remaining constituents of the
contraction \eqref{cont5} will be given in Complement II. As a
notational {\em road map\/} for the reader, we note at this stage
that Complement II involves an application of the ordinary
perturbation Lemma which will here yield, as an intermediate step,
yet another contraction of chain complexes, of the kind
\begin{equation*}
   \left(\Sigm^{\mathrm c}_{\delta}[sM]
     \begin{CD}
      \null @>{\nnabla}>> \null\\[-3.2ex]
      \null @<<{\widetilde \Pi}< \null
     \end{CD}
    \mathcal C[\fra g],\HH \right),
  \end{equation*}
to be given as \eqref{2.6a} below. In particular, $\delta$ is yet
another perturbation on $\Sigm^{\mathrm c}[sM]$ which we
distinguish in notation from the perturbation $\mathcal D$; apart
from trivial cases, the perturbation $\delta$ is {\em not\/}
compatible with the coalgebra structure on $\Sigm^{\mathrm
c}[sM]$, though, and the injection $\nnabla$ and homotopy $\HH$
differ from the ultimate injection $\overline \tau$ and homotopy
$H$.
\medskip

\noindent {\bf Complement I.} {\sl  The operator $\mathcal D$ and
twisting cochain $\tau$ are obtained as infinite series by the
following recursive procedure:
\begin{align}
\tau^1 &= \nabla \tau_{M}\colon\Sigm^{\mathrm c} \to \fra g,
\label{proc00}
\\
\tau^{\jj} &= \tfrac 12 h([\tau^1,\tau^{{\jj}-1}] +  \dots +
[\tau^{{\jj}-1},\tau^1])\colon\Sigm^{\mathrm c} \to \fra g,\  \jj
\geq 2, \label{proc0}
\\
\tau &= \tau^1 + \tau^2 + \dots \colon \Sigm^{\mathrm c} \to \fra
g,\label{proc1}
\\
\mathcal D &= \mathcal D^1 + \mathcal D^2 + \dots \colon
\Sigm^{\mathrm c} \to \Sigm^{\mathrm c} \label{proc2}
\end{align}
where, for $\jj \geq 1$, $\mathcal D^{\jj}$ is the coderivation of
$\Sigm^{\mathrm c}[sM]$ determined by the identity
\begin{equation}
\tau_{M} \mathcal D^{\jj} = \tfrac 12 \pi ([\tau^1,\tau^{\jj}] +
\dots + [\tau^{\jj},\tau^1]) \colon \Sigm_{\jj+ 1}^{\mathrm c} \to
M. \label{proc3}
\end{equation}
In particular, for $\jj \geq 1$, the coderivation $\mathcal
D^{\jj}$ is zero on $\mathrm F_j\Sigm^{\mathrm c}$ and lowers
coaugmentation filtration by $\jj$.}

\medskip

The sums \eqref{proc1} and \eqref{proc2} are in general infinite.
However, applied to a specific element which, since
$\Sigm^{\mathrm c}$ is cocomplete, necessarily lies in some finite
filtration degree subspace, since the operators $\mathcal D^{\jj}$
($\jj \geq 1$) lower coaugmentation filtration by $\jj$, only
finitely many terms will be non-zero, whence the convergence is
naive.

In the special case where the original contraction \eqref{cont1}
is the trivial contraction of the kind
\begin{equation}
   (\fra g
     \begin{CD}
      \null @>{\mathrm{Id}}>> \null\\[-3.2ex]
      \null @<<{\mathrm{Id}}< \null
     \end{CD}
    \fra g, 0) ,
\label{trivial}
  \end{equation}
$M$ and $\fra g$ coincide as chain complexes, the perturbation
$\ppartial$ coincides with the perturbation $\partial$ determined
by the graded Lie bracket $\ppppartial$ on $\fra g$, and
$\Sigm_{\ppartial}^{\mathrm c}[sM]$ coincides with the ordinary
{\scs CCE\/} or {\em classifying\/} coalgebra $\mathcal C[\fra g]$
for $\fra g$; the Lie algebra twisting cochain $\tau$ then comes
down to the {\scs CCE\/} or {\em universal\/} Lie algebra twisting
cochain $\tau_{\fra g}\colon \mathcal C[\fra g] \to \fra g$ for
$\fra g$ and in fact coincides with $\tau_1$ (in the present
special case) and, furthermore, the new contraction \eqref{cont5}
then amounts to the trivial contraction
\begin{equation*}
   (\mathcal
C [\fra g]
     \begin{CD}
      \null @>{\mathrm{Id}}>> \null\\[-3.2ex]
      \null @<<{\mathrm{Id}}< \null
     \end{CD}
    \mathcal
C [\fra g], 0) .
  \end{equation*}
In fact, in this case, the higher terms $\tau^{\jj}$ and
$\ppartial^{\jj}$ ($\jj \geq 2$) are obviously zero, and the
operator $\ppartial^1$ manifestly coincides with the CCE-operator.
Likewise, in the special case where the bracket on $\fra g$ is
trivial or, more generally, when $M$ carries a graded Lie bracket
in such a way that $\nabla$ is a morphism of differential graded
Lie algebras, the construction plainly stops after the first step,
and $\tau = \tau^1$.

\medskip

\noindent {\bf Complement II.} {\sl Application of the ordinary
perturbation Lemma (reproduced below as Lemma {\rm
\ref{ordinary}}) to the perturbation $\partial$ on $\Sigm^{\mathrm
c}[s\fra g]$ determined by the graded Lie algebra structure on
$\fra g$ and the induced {\em filtered\/} contraction
\begin{equation}
   \left(\Sigm^{\mathrm c}[sM]
     \begin{CD}
      \null @>{\Sigm^{\mathrm c}[s\nabla]}>> \null\\[-3.2ex]
      \null @<<{\Sigm^{\mathrm c}[s\pi]}< \null
     \end{CD}
    \Sigm^{\mathrm c}[s\fra g],\Sigm^{\mathrm c}[sh] \right)
\label{2.3}
  \end{equation}
of {\em coaugmented differential graded coalgebras\/}, the
filtrations being the ordinary coaugmentation filtrations, yields
the perturbation $\delta$ of the differential $d^0$ on
$\Sigm^{\mathrm c}[sM]$ and the contraction
\begin{equation}
   \left(\Sigm^{\mathrm c}_{\delta}[sM]
     \begin{CD}
      \null @>{\nnabla}>> \null\\[-3.2ex]
      \null @<<{\widetilde \Pi}< \null
     \end{CD}
    \mathcal C[\fra g],\HH \right)
\label{2.6a}
  \end{equation}
of chain complexes. Furthermore, the composite
\begin{equation}
\begin{CD}
\Phi\colon\Sigm^{\mathrm c}_{\pppartial}[sM] @>{\overline \tau}>>
\mathcal C[\fra g] @>{\widetilde \Pi}>> \Sigm^{\mathrm
c}_{\delta}[sM]
\end{CD}
\label{comp2}
\end{equation}
is an isomorphism of chain complexes, and the morphisms
\begin{align}
\Pi &= \Phi^{-1}\widetilde \Pi \colon \mathcal C[\fra g]
\longrightarrow \Sigm^{\mathrm c}_{\pppartial}[sM],
\label{Pipert1}
\\
H &= \HH- \HH \overline \tau\,\Pi \colon \mathcal C[\fra g]
\longrightarrow \mathcal C[\fra g] \label{Hpert1}
\end{align}
complete the construction of the contraction \eqref{cont5}.}

\medskip
In general, none of the morphisms $\delta$, $\nnabla$, $\widetilde
\Pi$, $\Pi$, $\HH$, $H$ is compatible with the coalgebra
structures. The isomorphism $\Phi$ admits an explicit description
in terms of the data as a {\em perturbation of the identity\/} and
so does its inverse; details will be given in Section
\ref{proofliepert} below.

\section{Some additional technical prerequisites}
\label{pre}

Let $Y$ be a chain complex. The {\em cofree\/} coaugmented
differential graded cocommutative coalgebra or {\em graded
symmetric\/} coalgebra $\Sigm^{\mathrm c}[Y]$ on the chain complex
$Y$ is characterized by a universal property as usual. To
guarantee the existence of a diagonal map for $\Sigm^{\mathrm
c}[Y]$, some hypothesis is necessary, though: The ordinary tensor
coalgebra $\mathrm T^{\mathrm c}[Y]$, that is, the cofree
(coaugmented) coalgebra on $Y$, decomposes as the direct sum
\[
\mathrm T^{\mathrm c}[Y] = \oplus_{j=0}^{\infty} \mathrm
T_j^{\mathrm c}[Y]
\]
of its homogeneous constituents $\mathrm T_j^{\mathrm c}[Y] =
Y^{\otimes j}$ ($j \geq 0$). For $j \geq 0$, let $\Sigm_j^{\mathrm
c}[Y]\subseteq \mathrm T_j^{\mathrm c}[Y]$ be the submodule of
invariants in the $j$'th tensor power $\mathrm T_j^{\mathrm c}[Y]$
relative to the obvious action on $\mathrm T_j^{\mathrm c}[Y]$ of
the symmetric group $S_j$ on $j$ letters, and let $\Sigm^{\mathrm
c}[Y]$ be the direct sum
\[ \Sigm^{\mathrm
c}[Y]=\oplus_{j=0}^{\infty} \Sigm_j^{\mathrm c}[Y]
\]
of chain complexes. So far, the construction is completely
general, even functorial, and works over any ground ring. In
particular, a chain map $\phi \colon Y_1 \to Y_2$ induces a chain
map
\[
\Sigm^{\mathrm c}[\phi] \colon\Sigm^{\mathrm c}[Y_1]
\longrightarrow \Sigm^{\mathrm c}[Y_2].
\]
However, some hypothesis is, in general, necessary in order for
the homogeneous constituents
\[
\mathrm T_{j+k}^{\mathrm c}[Y] \longrightarrow \mathrm
T_j^{\mathrm c}[Y] \otimes \mathrm T_k^{\mathrm c}[Y]\ (j,k \geq
0)
\]
of the diagonal map $\Delta \colon \mathrm T^{\mathrm c}[Y] \to
\mathrm T^{\mathrm c}[Y] \otimes \mathrm T^{\mathrm c}[Y]$ of the
graded tensor coalgebra $\mathrm T^{\mathrm c}[Y]$ to induce a
graded diagonal map on $\Sigm^{\mathrm c}[Y]$. (I am indebted to
P. Urbanski for having prodded me to clarify this point.) Indeed,
the diagonal map for $Y$ induces a morphism from $\Sigm^{\mathrm
c}[Y]$ to $\Sigm^{\mathrm c}[Y\oplus Y]$ and the diagonal map for
$\Sigm^{\mathrm c}[Y]$ is well defined whenever the canonical
morphism
\begin{equation}
\oplus _{j_1 + j_2 = k}\Sigm_{j_1}^{\mathrm c}[Y] \otimes
\Sigm_{j_2}^{\mathrm c}[Y] \longrightarrow \Sigm_k^{\mathrm
c}[Y\oplus Y] \label{canm} \end{equation}
is an isomorphism for
every $k\geq 1$.

To explain the basic difficulty, let $k \geq 1$, let $Y_1$ and
$Y_2$ be two graded $R$-modules, and consider the $k$'th
homogeneous constituent
\[
\Sigm_{k}^{\mathrm c}[Y_1 \oplus Y_2] = \left(\left(Y_1 \oplus
Y_2\right)^{\otimes k}\right)^{S_k} \subseteq \left(Y_1 \oplus
Y_2\right)^{\otimes k}
\]
of $\Sigm^{\mathrm c}[Y_1 \oplus Y_2]$. For $ 0 \leq j \leq k$,
let
\[
\Sigm_{\binom kj}[Y_1 \oplus Y_2] =Y_1^{\otimes j} \otimes
Y_2^{\otimes (k-j)} \oplus Y_1^{\otimes (j-1)} \otimes Y_2\otimes
Y_1\otimes Y_2^{\otimes (k-j-1)} \oplus \ldots \oplus Y_1^{\otimes
(k-j)} \otimes Y_2^{\otimes j},
\]
the direct sum of $\binom kj$ summands which arises by
substituting in the possible choices of $j$ objects out of $k$
objects a tensor factor of $Y_1$ for each object and filling in
the \lq\lq holes\rq\rq\ remaining between the various tensor
powers of $Y_1$ by the appropriate tensor powers of $Y_2$. Let
$RS_k$ denote the group ring of $S_k$. For $ 0 \leq j \leq k$,
relative to the $S_j$-and $S_{k-j}$-actions on $Y_1^{\otimes j}$
and $Y_2^{\otimes (k-j)}$, respectively, there is a canonical
isomorphism
\[
\left(Y_1^{\otimes j} \otimes Y_2^{\otimes (k-j)}\right)
\otimes_{S_j \times S_{k-j}}RS_k \longrightarrow \Sigm_{\binom
kj}[Y_1 \oplus Y_2]
\]
of $S_k$-modules. Consequently, for $ 0 \leq j \leq k$, the
canonical injection
\[
\left(Y_1^{\otimes j} \otimes Y_2^{\otimes (k-j)}\right)^{S_j
\times S_{k-j}}\longrightarrow \left(\Sigm_{\binom kj}[Y_1 \oplus
Y_2]\right)^{S_k}
\]
is an isomorphism. As an $S_k$-module, $\left(Y_1 \oplus
Y_2\right)^{\otimes k}$ is well known to decompose as the direct
sum
\[
\left(Y_1 \oplus Y_2\right)^{\otimes k}
=\oplus_{j=0}^k\Sigm_{\binom kj}[Y_1 \oplus Y_2]
\]
whence $\Sigm_k^{\mathrm c}[Y_1 \oplus Y_2]$ decomposes as the
direct sum
\[
\Sigm_k^{\mathrm c}[Y_1 \oplus Y_2] = \oplus_{j=0}^k
\left(Y_1^{\otimes j} \otimes Y_2^{\otimes (k-j)}\right)^{S_j
\times S_{k-j}} .
\]
However, some hypothesis is needed in order for the canonical
morphisms
\[
\Sigm_j^{\mathrm c}[Y_1] \otimes \Sigm_{k-j}^{\mathrm c}[Y_2] =
\left(Y_1^{\otimes j}\right)^{S_j} \otimes \left(Y_2^{\otimes
(k-j)}\right)^{S_{k-j}} \longrightarrow \left(Y_1^{\otimes j}
\otimes Y_2^{\otimes (k-j)}\right)^{S_j \times S_{k-j}}
\]
of graded $R$-modules  to be  isomorphisms for $1 \leq j \leq
k-1$.

To explain an important special case where the cofree graded
cocommutative coalgebra exists we recall that, when
$\Sigm^{\mathrm c}[Y]$ exists, the addition map of $Y$ induces a
multiplication map for $\Sigm^{\mathrm c}[Y]$ turning
$\Sigm^{\mathrm c}[Y]$ into a differential graded Hopf algebra. In
particular, when $Y$ is free and concentrated in even degrees,
$\Sigm^{\mathrm c}[Y]$ amounts to the familiar {\em divided powers
Hopf algebra\/} $\Gamma[Y]$ on a basis of $Y$, cf.
\cite{cartanse}; the underlying graded $R$-module is then free as
well, that is, $\Sigm^{\mathrm c}[Y]$ is then free as a graded
$R$-module. This kind of construction extends to graded projective
modules, that is, $\Sigm^{\mathrm c}[Y]$ exists as a coaugmented
graded cocommutative coalgebra whenever $Y$, concentrated in even
degrees, is projective as a graded $R$-module, and each
homogeneous constituent $\Sigm_k^{\mathrm c}[Y]$ ($k \geq 1$) is
then projective as well. Given a general graded module $Y$ which
is projective as a graded $R$-module, denote the even and odd
constituents by $Y_{\mathrm{ev}}$ and $Y_{\mathrm{odd}}$,
respectively; the coaugmented graded cocommutative coalgebra
$\Sigm^{\mathrm c}[Y]$ on $Y$ can then be constructed as the
tensor product
\[
\Sigm^{\mathrm c}[Y] = \Sigm^{\mathrm c}[Y_{\mathrm{ev}}] \otimes
\Sigm^{\mathrm c}[Y_{\mathrm{odd}}]
\]
where $\Sigm^{\mathrm c}[Y_{\mathrm{odd}}]$ is the coalgebra which
underlies the ordinary exterior Hopf algebra on
$Y_{\mathrm{odd}}$. Thus, to sum up, {\sl given a chain complex
$Y$ which, as a graded $R$-module, is projective, the coaugmented
differential graded cocommutative coalgebra $\Sigm^{\mathrm c}[Y]$
on $Y$ exists\/}.

Likewise, given a general chain complex $Y$, under suitable
circumstances, the diagonal map $Y \to Y \oplus Y \cong Y \times
Y$ of $Y$ induces a diagonal map for the graded symmetric algebra
$\Sigm[Y]$ on $Y$ turning $\Sigm[Y]$ into a differential graded
Hopf algebra and, whenever $\Sigm^{\mathrm c}[Y]$ exists, the
identity morphism of $Y$ induces a morphism
$\Sigm[Y]\to\Sigm^{\mathrm c}[Y]$ of differential graded Hopf
algebras. Over a field of characteristic zero, this morphism is an
isomorphism (since, on the homogeneous degree $j$ constituents
($j\geq 1$), division by the factorial $j!$ is possible) and, in
the literature, over a field of characteristic zero, the coalgebra
underlying the graded symmetric Hopf algebra $\Sigm[Y]$ is usually
taken as a model for the coaugmented differential graded
cocommutative coalgebra on $Y$. More generally, suppose that the
ground ring $R$ contains the rational numbers as a subring. Then
the coalgebra which underlies $\Sigm[Y]$ may be taken as a model
for $\Sigm^{\mathrm c}[Y]$, cf. the discussion in Appendix B of
\cite{quilltwo}. Indeed, the map
\[
\Sigm[Y] \longrightarrow \mathrm T ^{\mathrm c}[Y],\ x_1\ldots x_n
\longmapsto \frac 1{n!}\sum _{\sigma\in S_n} x_{\sigma 1} \otimes
\ldots \otimes x_{\sigma n},\ x_j \in Y,\ n>0,
\]
induces an explicit isomorphism of $\Sigm[Y]$ onto $\Sigm^{\mathrm
c}[Y]\subseteq \mathrm T ^{\mathrm c}[Y]$. However, over the
integers, given a free graded abelian group $Y$ concentrated in
even degrees, the resulting morphism $\Sigm[Y]\to\Gamma[Y]$ of
graded Hopf algebras is {\em not\/} an isomorphism. Thus the above
construction of the coaugmented differential graded cocommutative
coalgebra $\Sigm^{\mathrm c}[Y]$ on $Y$ is more general.

As {\em chain complexes\/}, not just as graded objects, the
Hom-complexes $\mathrm{Hom}(\Sigm^{\mathrm c},\fra g)$ and
$\mathrm{Hom}(\mathrm F_n\Sigm^{\mathrm c},\fra g)$ ($n\geq 0$)
manifestly decompose as direct products
\[
\mathrm{Hom}(\Sigm^{\mathrm c},\fra g)\cong \prod_{j=0}^{\infty}
\mathrm{Hom}(\Sigm_j^{\mathrm c},\fra g),\quad
\mathrm{Hom}(\mathrm F_n\Sigm^{\mathrm c},\fra g)\cong
\prod_{j=0}^n \mathrm{Hom}(\Sigm_j^{\mathrm c},\fra g),\ (n\geq
0),
\]
and restriction mappings induce a sequence
\begin{equation}
\ldots \longrightarrow \mathrm{Hom}(\mathrm F_{n+1}\Sigm^{\mathrm
c},\fra g) \longrightarrow \mathrm{Hom}(\mathrm F_n\Sigm^{\mathrm
c},\fra g)\longrightarrow \ldots \longrightarrow
\mathrm{Hom}(\mathrm F_1\Sigm^{\mathrm c},\fra g)\longrightarrow
\fra g \label{proj1}
\end{equation}
of surjective morphisms of differential graded Lie algebras.
Furthermoire, by construction, for each $n\geq 0$, the canonical
injection of $\mathrm F_n\Sigm^{\mathrm c}$ into $\Sigm^{\mathrm
c}$ is a morphism of coaugmented differential graded coalgebras
and hence induces a projection $\mathrm{Hom}(\Sigm^{\mathrm
c},\fra g) \to \mathrm{Hom}(\mathrm F_n\Sigm^{\mathrm c},\fra g)$
of differential graded Lie algebras, and these projections
assemble to an isomorphism from $\mathrm{Hom}(\Sigm^{\mathrm
c},\fra g)$ onto the projective limit $\varprojlim
\mathrm{Hom}(\mathrm F_n\Sigm^{\mathrm c},\fra g)$ of
\eqref{proj1} in such a way that, in each degree, the limit is
attained at a finite stage.

\section{The crucial step}

For ease of exposition, we introduce the notation
$\pppartial^0=\pppartial_0=0$. For $a\geq 1$, let
\begin{align}
\tau_a&= \tau^1 + \tau^2 + \dots + \tau^a,
\\
\ppartial_a&= \ppartial^1 + \ppartial^2 + \dots + \ppartial^a,
\\
\Theta_{a+1}&=-D\tau_a - \tau_a\pppartial_{a-1}+ \tfrac12
[\tau_a,\tau_a]=-d\tau_a - \tau_a (d^0+\pppartial_{a-1}) +
\tfrac12 [\tau_a,\tau_a] \colon \Sigm^{\mathrm c} \longrightarrow
\fra g; \label{twist111}
\\
\vartheta_{a+1}&= \Theta_{a+1}|_{\Sigm_{a+1}^{\mathrm c}}\colon
\Sigm_{a+1}^{\mathrm c} \longrightarrow \fra g. \label{twist444}
\end{align}

The crucial step for the proof of the Lie algebra perturbation
Lemma, in particular for the statement given as Complement I
above, is provided by the following.

\begin{Lemma} \label{lem12} Let $a \geq 1$.
\begin{align}
\pi\tau^{a+1}&=0; \label{twist3}
\\
\vartheta_{a+1}&= -\tau^2 \ppartial^{a-1}-\ldots -\tau^a
\ppartial^{1}+\tfrac 12\left([\tau^1,\tau^a]+ \ldots +
[\tau^a,\tau^1]\right)\colon \Sigm_{a+1}^{\mathrm c}
\longrightarrow \fra g;
\label{twist331}\\
h\vartheta_{a+1}&=\tfrac 12 h\left([\tau^1,\tau^a]+ \ldots +
[\tau^a,\tau^1]\right) = \tau^{a+1}\colon \Sigm_{a+1}^{\mathrm c}
\longrightarrow \fra g; \label{twist332}
\\
\pi\vartheta_{a+1}&=\tfrac 12 \pi\left([\tau^1,\tau^a]+ \ldots +
[\tau^a,\tau^1]\right) = \tau_M\ppartial^a\colon
\Sigm_{a+1}^{\mathrm c} \longrightarrow M ;
\label{twist333}
\\
\Theta_{a+1}& \colon \Sigm^{\mathrm c} \longrightarrow \fra g
\text{\ is zero on\ } \mathrm F_a\Sigm^{\mathrm c}, \text{\ i.~e.
goes to zero in\ }\mathrm{Hom}(\mathrm F_a \Sigm^{\mathrm c},\fra
g); \label{twist1}
\\
D\vartheta_{a+1}& =\tau_1(\ppartial^1\ppartial^{a-1} + \dots +
\ppartial^{a-1}\ppartial^1) \colon \Sigm_{a+1}^{\mathrm c}
\longrightarrow \fra g; \label{twist4}
\\
D\tau^{a+1}&=\vartheta_{a+1} - \tau^1 \ppartial^a \colon
\Sigm^{\mathrm c} \longrightarrow \fra g; \label{twist5}
\\
0&=\dd^0 \ppartial^{a} + \ppartial^1\ppartial^{a-1}+ \dots +
\ppartial^{a-1}\ppartial^1+ \ppartial^{a}\dd^0 \colon
\Sigm^{\mathrm c} \longrightarrow \Sigm^{\mathrm c}.
\label{twist22}
\end{align}
\end{Lemma}

For clarity we note that, for $a=1$, the formula \eqref{twist331}
signifies
\begin{equation} \vartheta_2=\tfrac 12[\tau^1,\tau^1],
\label{clar}
\end{equation}
and \eqref{twist4} signifies that $\vartheta_2$ is a cycle, i.~e.
$D \vartheta_2 = 0$.

Let $b \geq 1$. The property \eqref{twist22} implies that, on the
differential graded coalgebra $\mathrm F_{b+1}\Sigm^{\mathrm c}$,
the operator $\pppartial_{b-1}$ is a coalgebra perturbation of the
differential $d^0$ and hence $d^0+\pppartial_{b-1}$ is a coalgebra
differential and thence turns the coaugmented graded cocommutative
coalgebra $\mathrm F_{b+1}\Sigm^{\mathrm c}$ into a coaugmented
differential graded cocommutative coalgebra which, according to
the convention introduced above, we write as $\mathrm
F_{b+1}\Sigm^{\mathrm c}_{\pppartial_{b-1}}$. Furthermore,
property \eqref{twist1} implies that the restriction of $\tau_b$
to $\mathrm F_b\Sigm^{\mathrm c}$ is a Lie algebra twisting
cochain $\mathrm F_b\Sigm^{\mathrm c}_{\pppartial_{b-1}} \to \fra
g$ whence, as $b$ tends to infinity, $\tau_b$ tends to a Lie
algebra twisting cochain, that is, $\tau$ is a Lie algebra
twisting cochain. Indeed, in a given degree, the statements of
Lemma \ref{lem12} come down to corresponding statements in a
suitable finite stage constituent of the sequence \eqref{proj1}.

\begin{proof} The property \eqref{twist3} is an immediate consequence of the
annihilation properties \eqref{side}. Next, let $a \geq 1$.  For
degree reasons, the restriction of
\[
\Theta_{a+1}=-D\tau_a - \tau_a\pppartial_{a-1}+ \tfrac12
[\tau_a,\tau_a]=-d\tau_a - \tau_a (d^0+\pppartial_{a-1}) +
\tfrac12 [\tau_a,\tau_a] \colon \Sigm^{\mathrm c} \longrightarrow
\fra g
\]
to $\Sigm_{a+1}^{\mathrm c}$ comes down to
\[
-\tau^2 \ppartial^{a-1}-\ldots -\tau^a \ppartial^{1}+\tfrac
12\left([\tau^1,\tau^a]+ \ldots + [\tau^a,\tau^1]\right) \colon
\Sigm_{a+1}^{\mathrm c} \longrightarrow \fra g,
\]
whence \eqref{twist331}, being interpreted as \eqref{clar} for
$a=1$. The identity \eqref{twist331}, combined with the
annihilation properties \eqref{side}, immediately implies
 \eqref{twist332} and
\eqref{twist333}, in view of the definitions \eqref{proc1} and
\eqref{proc3} of the terms $\tau^{j+1}$ and $\ppartial^j$,
respectively, for $j \geq 1$.

Furthermore, the property \eqref{twist5} is a formal consequence
of the definitions \eqref{proc3} and \eqref{twist444}, combined
with \eqref{twist4} and the annihilation properties \eqref{side}.
Indeed,
\begin{align*}
D\tau^{a+1}&= D(h\vartheta_{a+1}) = (Dh)\vartheta_{a+1}
-hD\vartheta_{a+1}
\\&=
\vartheta_{a+1} - \nabla \pi \vartheta_{a+1} -h \tau_1
(\ppartial^1\ppartial^{a-1} + \dots + \ppartial^{a-1}\ppartial^1)
\\
&= \vartheta_{a+1} - \nabla \tau_M \ppartial^a
\\
&= \vartheta_{a+1} - \tau_1 \ppartial^a \colon
\Sigm_{a+1}^{\mathrm c} \longrightarrow \fra g,
\end{align*}
whence \eqref{twist5}.

By induction on $a$, we now establish the remaining assertions
\eqref{twist1}, \eqref{twist4}, and \eqref{twist22}. To begin
with, let $a=1$. Since $\tau_1$ is a cycle in
$\mathrm{Hom}(\Sigm^{\mathrm c},\fra g)$ and since
$[\tau_1,\tau_1]$ vanishes on $\mathrm F_1\Sigm^{\mathrm c}$,
\[ \Theta_2 =-d\tau_1 -\tau_1 \dd^0
+\tfrac 12[\tau_1,\tau_1] = \tfrac 12[\tau_1,\tau_1]\colon
\Sigm^{\mathrm c} \to \fra g
\]
vanishes on $\mathrm F_1\Sigm^{\mathrm c}$, whence \eqref{twist1}
holds for $a=1$. Furthermore, since $\Theta_2$ is a cycle, so is
$\vartheta_2$ whence \eqref{twist4} is satisfied for $a=1$.
Consequently
\begin{align*}
D\tau^2 &= Dh \vartheta_2 = \vartheta_2 - \nabla\pi\vartheta_2 =
\tfrac 12[\tau_1,\tau_1] - \nabla \pi \vartheta_2
\\
&= \tfrac 12[\tau_1,\tau_1] - \nabla \tau_M \ppartial^1 =\tfrac
12[\tau_1,\tau_1] -\tau_1 \ppartial^1
\end{align*}
whence \eqref{twist5} for $a=1$. Finally, the identity
\eqref{twist22} for $a=1$ reads
\begin{equation}
\dd^0 \ppartial^1 + \ppartial^1 \dd^0 = 0. \label{pert22}
\end{equation}
The identity \eqref{pert22}, in turn, is a consequence of the
bracket on $\fra g$ being compatible with the differential $d$ on
$\fra g$ since this compatibility entails that $[\tau_1,\tau_1]$
is a cycle in $\mathrm{Hom}(\Sigm^{\mathrm c},\fra g)$. Indeed,
since
\[
\dd^0 \ppartial^1 + \ppartial^1 \dd^0 =[\dd^0, \ppartial^1]
\]
is a coderivation on $\Sigm^{\mathrm c}$, the bracket being the
commutator bracket in the {\em graded Lie algebra of
coderivations\/} of $\Sigm^{\mathrm c}$, it suffices to show that
\[
\tau_M(\dd^0 \ppartial^1 + \ppartial^1 \dd^0) =\tau_M [\dd^0,
\ppartial^1]
\]
vanishes. However, since $\tau_M \ppartial^1 = \tfrac 12 \pi
[\tau_1,\tau_1]$, cf. \eqref{proc3},
\[
\tau_M(\dd^0 \ppartial^1 + \ppartial^1 \dd^0) = -d \tau_M
\ppartial^1 + \tau_M \ppartial^1 \dd^0 =-\tfrac 12 D (\pi
[\tau_1,\tau_1]),
\]
and $D (\pi [\tau_1,\tau_1])$ vanishes since $[\tau_1,\tau_1]$ is
a cycle in $\mathrm{Hom}(\Sigm^{\mathrm c},\fra g)$ and since
$\pi$ is a chain map. Consequently the identity \eqref{twist22}
holds for $a=1$. Thus the induction starts.

Even though this is not strictly necessary we now explain the case
$a=2$. This case is particularly instructive. Now
\begin{align*}
d\tau_2 + \tau_2 (d^0+ \pppartial^1) &= d\tau^1+d\tau^2 + \tau^1
\dd^0 +\tau^1 \ppartial^1 +\tau^2\dd^0 +\tau^2\ppartial^1
\\
&=D\tau_1 +D\tau^2 + \tau_1 \ppartial^1 +\tau^2\ppartial^1
\\
&=\tfrac 12[\tau_1,\tau_1]  - \tau_1\ppartial^1 + \tau_1
\ppartial^1 +\tau^2\ppartial^1
\\
&=\tfrac 12[\tau_1,\tau_1] +\tau^2\ppartial^1
\end{align*}
whence
\[
\Theta_3 =-d\tau_2 - \tau_2 (d^0+\pppartial^1) +\tfrac
12[\tau_2,\tau_2] = -\tau^2\ppartial^1 +[\tau^1,\tau^2] +\tfrac
12[\tau^2,\tau^2]
\]
which clearly vanishes on $\mathrm F_2\Sigm^{\mathrm c}$, whence
\eqref{twist1} holds for $a=2$. Furthermore, it is manifest that
the restriction $\vartheta_3$ of $\Theta_3$ to $\Sigm_3^{\mathrm
c}$ takes the form
\begin{equation}
\vartheta_3 =  [\tau^1,\tau^2]- \tau^2 \ppartial^1\colon
\Sigm_3^{\mathrm c} \longrightarrow \fra g, \label{vartheta3}
\end{equation}
which amounts to \eqref{twist331} for the special case $a=2$.
Hence
\begin{align*}
D\vartheta_3 &= -[\tau^1,D\tau^2]- (D\tau^2) \ppartial^1
\\
&= -[\tau^1,-\tau^1\ppartial^1 + \tfrac 12
[\tau^1,\tau^1]]+\tau^1\ppartial^1\ppartial^1 - \tfrac 12
[\tau^1,\tau^1]\ppartial^1
\\
&= [\tau^1,\tau^1\ppartial^1] +\tau^1\ppartial^1\ppartial^1 -
[\tau^1,\tau^1\ppartial^1]
\\
&=\tau^1\ppartial^1\ppartial^1
\end{align*}
whence \eqref{twist4} at stage $a=2$. Since
\[ D\tau^2 =\tfrac 12[\tau_1,\tau_1] -\tau_1 \ppartial^1
\]
and since $[\tau_1,[\tau_1,\tau_1]]=0$,
\begin{equation*}
D\left[\tau_1, \tau^2\right] = \left[\tau_1, -D\tau^2\right] =
\left[\tau_1,\tau_1 \ppartial_1 -\tfrac 12
\left[\tau_1,\tau_1\right]\right] = \left[\tau_1,\tau_1
\ppartial_1\right].
\end{equation*}
Thus, in view of \eqref{proc3}, viz.
\[
\pi\tau_1\ppartial_1 = \tfrac 12 \pi[\tau_1,\tau_1],
\]
we find
\begin{equation}
D(\pi\left [\tau_1,\tau^2 \right ]) = \pi \tau_1 \ppartial^1
\ppartial^1 , \label{del129}
\end{equation}
whence
\begin{equation}
D(\pi\left [\tau_1,\tau^2 \right ]) = \tau_M \ppartial^1
\ppartial^1 . \label{del130}
\end{equation}

Since, in view of \eqref{vartheta3} or \eqref{twist331},
\[
\pi \vartheta_3 = \pi[\tau_1,\tau^2] \colon \Sigm_3^{\mathrm
c}[sM] \longrightarrow M,
\]
$\ppartial^2 \colon \Sigm^{\mathrm c} \to \Sigm^{\mathrm c}$ is
the coderivation which is determined by the requirement that the
identity
\begin{equation}
\pi \vartheta_3 = \tau_M \ppartial^2\colon \Sigm_3^{\mathrm c}[sM]
\longrightarrow M\label{del1011}
\end{equation}
be satisfied. Then
\begin{equation}
\tau_M (\dd^0\ppartial^2 + \ppartial^1\ppartial^1 + \ppartial^2
\dd^0 )=\tau_M (D\ppartial^2 + \ppartial^1\ppartial^1)=0 .
\end{equation}
Indeed,
\begin{align*}
\tau_M D\ppartial^2&=-D(\tau_M \ppartial^2) \\
&=-D(\pi\vartheta_3) =-D(\pi\left [\tau_1,\tau^2 \right ]) \\
&= -\tau_M \ppartial^1 \ppartial^1 .
\end{align*}
Consequently
\begin{equation}
\dd^0\ppartial^2 + \ppartial^1\ppartial^1 + \ppartial^2 \dd^0=0
\label{asso2}
\end{equation}
since $\dd^0\ppartial^2 + \ppartial^1\ppartial^1 + \ppartial^2
\dd^0$ is a coderivation of $\Sigm^{\mathrm c}[sM]$. This
establishes the identity \eqref{twist22} for $a=2$.

We pause for the moment; suppose that we are in the special
situation where the original contraction \eqref{cont1} is the
trivial one of the kind \eqref{trivial} and identify $M$ with the
chain complex which underlies $\fra g$, endowed with the zero
bracket. Then $\Sigm^{\mathrm c}$ amounts to the {\scs
CCE\/}-coalgebra for $M$ (endowed with the zero bracket), the
operator $\ppartial^1$ is precisely the ordinary {\scs
CCE\/}-operator relative to the Lie bracket on $\fra g$, the
twisting cochain $\tau_1$ is the {\scs CCE\/}-twisting cochain
relative to the Lie bracket on $\fra g$, the term $\tau^2$ is zero
(since $h$ is zero), and the construction we are in the process of
explaining stops at the present stage. Indeed, $\tau_1$ then
coincides with $\tau_2$ and $\Theta_3=0$. Moreover, the identity
\begin{equation}
\ppartial^1 \ppartial^1 =0  \label{pert3}
\end{equation}
is then equivalent to the bracket on $\fra g$ satisfying the
graded Jacobi identity.

Likewise, in the special case where the differential on $M$ is
zero so that $M$ amounts to the homology $\mathrm H(\fra g)$ of
$\fra g$, the identity \eqref{asso2} comes down to
\begin{equation}
\ppartial^1\ppartial^1 =0. \label{asso22}
\end{equation}
This identity, in turn, is then equivalent to the fact that the
induced graded bracket on $\mathrm H(\fra g)$ satisfies the graded
Jacobi identity.

We now return to the case of a general contraction \eqref{cont1}.
Let $b> 2$ and suppose, by induction that, at stage $a$, $2 \leq a
< b$, \eqref{twist1} -- \eqref{twist22} have been established. Our
aim is to show that \eqref{twist1} -- \eqref{twist22} hold at
stage $b$. Now
\begin{align*}
\Theta_{b+1} &=-D \tau_b -\tau_b \pppartial_{b-1} + \tfrac
12[\tau_b,\tau_b] \colon \Sigm^{\mathrm c}[sM] \longrightarrow
\fra g\\
&= -D \tau_{b-1} -D\tau^b -(\tau_{b-1} + \tau^b) (\pppartial_{b-2}
+ \ppartial^{b-1})
\\&\qquad
 + \tfrac 12([\tau_{b-1},\tau_{b-1}]
+[\tau^b,\tau_{b-1}]+[\tau_{b-1},\tau^b]+[\tau^b,\tau^b])
\\
&= \Theta_b - D\tau^b -\tau_{b-1} \pppartial^{b-1} -\tau^b
\pppartial_{b-1} + \tfrac
12([\tau^b,\tau_{b-1}]+[\tau_{b-1},\tau^b]+[\tau^b,\tau^b]).
\end{align*}
By the inductive hypothesis \eqref{twist1} at stage $b-1$,
$\Theta_b$ vanishes on $\mathrm F_{b-1} \Sigm^{\mathrm c}$ whence
$\Theta_{b+1}$ vanishes on $\mathrm F_{b-1} \Sigm^{\mathrm c}$ as
well since the remaining terms obviously vanish on $\mathrm
F_{b-1} \Sigm^{\mathrm c}$. Moreover,
\[
\Theta_{b+1}\big|_{\Sigm_b^{\mathrm c}}= \vartheta_b - D\tau^b
-\tau_{b-1} \pppartial^{b-1} -\tau^b \pppartial_{b-1} + \tfrac
12([\tau^b,\tau_{b-1}]+[\tau_{b-1},\tau^b]+[\tau^b,\tau^b]).
\]
In view of the inductive hypothesis \eqref{twist5},
\[
D\tau^b=\vartheta_b - \tau^1 \ppartial^{b-1}
\]
and, for degree reasons,
\[
\tau^1 \ppartial^{b-1} = \tau_{b-1} \ppartial^{b-1}
\]
whence
\[
\Theta_{b+1}\big|_{\Sigm_b^{\mathrm c}}=  -\tau^b \pppartial_{b-1}
+ \tfrac
12([\tau^b,\tau_{b-1}]+[\tau_{b-1},\tau^b]+[\tau^b,\tau^b])
\]
which, for degree reasons, is manifestly zero. Consequently
$\Theta_{b+1}$ vanishes on $\mathrm F_b \Sigm^{\mathrm c}$ whence
\eqref{twist1} at stage $b$.

Next we establish the identity \eqref{twist4} at stage $b$. Recall
that, by construction, cf. \eqref{twist111},
\[
\Theta_{b+1} = -D\tau_b - \tau_b\pppartial_{b-1} + \tfrac 12
[\tau_b,\tau_b]
\]
whence
\[
D\Theta_{b+1} = -D(\tau_b\pppartial_{b-1}) + \tfrac 12
D[\tau_b,\tau_b] = -(D\tau_b)\pppartial_{b-1}
+\tau_bD\pppartial_{b-1} + [D\tau_b,\tau_b].
\]
However, $b\geq 2$  whence $\pppartial_{b-1}$ lowers filtration,
i.~e. maps $\mathrm F_{b+1}\Sigm^{\mathrm c}$ to $\mathrm
F_b\Sigm^{\mathrm c}$. Since\eqref{twist1} has already been
established at stage $b$, restricted to $\mathrm
F_{b+1}\Sigm^{\mathrm c}$,
\[
0=-\Theta_{b+1}\ppartial_{b-1} =D\tau_b\ppartial_{b-1}
+\tau_b\pppartial_{b-1}\ppartial_{b-1} - \tfrac
12[\tau_b,\tau_b]\ppartial_{b-1}
\]
whence
\[
(D\tau_b)\pppartial_{b-1}+ \tau_b\pppartial_{b-1}\pppartial_{b-1}
= \tfrac 12[\tau_b,\tau_b] \pppartial_{b-1}= [\tau_b,\tau_b
\pppartial_{b-1}]\in \mathrm{Hom}(\mathrm F_{b+1}\Sigm^{\mathrm
c},\fra g).
\]
Consequently
\begin{align*}
D\vartheta_{b+1} &= \tau_b\pppartial_{b-1}\pppartial_{b-1}
-[\tau_b,\tau_b\pppartial_{b-1}] +\tau_bD\pppartial_{b-1} +
[D\tau_b,\tau_b] \in \mathrm{Hom}(\mathrm F_{b+1}\Sigm^{\mathrm
c},\fra g)
\\
&= \tau_b(\pppartial_{b-1}\pppartial_{b-1}+D\pppartial_{b-1}) -
[\tau_b,\tau_b\pppartial_{b-1}+D\tau_b].
\end{align*}
By induction, in view of \eqref{twist1},
\[
\tau_b\pppartial_{b-1}+ D\tau_b=\tfrac 12[\tau_b,\tau_b] \in
\mathrm{Hom}(\mathrm F_b\Sigm^{\mathrm c},\fra g),
\]
whence
\[
[\tau_b,\tau_b\pppartial_{b-1}+D\tau_b]= \tfrac
12[\tau_b,[\tau_b,\tau_b]]
\]
which, each homogeneous constituent of $\tau_b$ being odd, is
zero, in view of the graded Jacobi identity in
$\mathrm{Hom}(\mathrm F_b\Sigm^{\mathrm c},\fra g)$. Moreover, by
induction, in view of \eqref{twist22}, for $1 \leq a < b$,
\[
\ppartial^1\ppartial^{a-1}+ \dots + \ppartial^{a-1}\ppartial^1+
D\ppartial^{a} = 0
\]
whence
\[
\pppartial_{b-1}\pppartial_{b-1}+D\pppartial_{b-1}=
\ppartial^1\ppartial^{b-1} + \dots + \ppartial^{b-1}\ppartial^1.
\]
Consequently, on $\Sigm_{b+1}^{\mathrm c}$,
\[
D\vartheta_{b+1} = \tau_b(\ppartial^1\ppartial^{b-1} + \dots +
\ppartial^{b-1}\ppartial^1)\colon \Sigm_{b+1}^{\mathrm c}\to \fra
g
\]
and, for degree reasons,
\[
\tau_b(\ppartial^1\ppartial^{b-1} + \dots +
\ppartial^{b-1}\ppartial^1) = \tau_1(\ppartial^1\ppartial^{b-1} +
\dots + \ppartial^{b-1}\ppartial^1)\colon \Sigm_{b+1}^{\mathrm
c}\to \fra g
\]
whence
\[
D\vartheta_{b+1} = \tau_1(\ppartial^1\ppartial^{b-1} + \dots +
\ppartial^{b-1}\ppartial^1).
\]
This establishes the identity \eqref{twist4} at stage $b$.

Alternatively, in view of what has already been proved, by virtue
of \eqref{twist331},
\[
\vartheta_{b+1}= -\tau^2 \ppartial^{b-1}-\ldots -\tau^b
\ppartial^{1}+\tfrac 12\left([\tau^1,\tau^b]+ \ldots +
[\tau^b,\tau^1]\right)\colon \Sigm_{b+1}^{\mathrm c}
\longrightarrow \fra g,
\]
whence, since $D\tau^1=0$,
\begin{equation}
\begin{aligned}
D\vartheta_{b+1}&= -(D\tau^2) \ppartial^{b-1} + \tau^2
D\ppartial^{b-1} \pm \ldots -(D\tau^{b-1}) \ppartial^2
\\
&\quad
 +\tau^{b-1}
D\ppartial^2 -D\tau^b \ppartial^1
\\
&\quad +[D\tau^2,\tau^{b-1}]+ \ldots + [D\tau^b,\tau^1] .
\end{aligned}
\label{alt}
\end{equation}
Thus, using the inductive hypotheses, we can establish
\eqref{twist4} at stage $b$ by evaluating the terms on the
right-hand side of \eqref{alt}.

Finally, to settle the identity \eqref{twist22} at stage $b$, we
note first that, since
\[
\dd^0 \ppartial^{b} + \ppartial^1\ppartial^{b-1}+ \dots +
\ppartial^{b-1}\ppartial^1+ \ppartial^{b}\dd^0 =
\ppartial^1\ppartial^{b-1}+ \dots + \ppartial^{b-1}\ppartial^1+
D\ppartial^{b}
\]
is a coderivation of $\Sigm^{\mathrm c}$, it suffices to prove
that
\[
\tau_M\left(\ppartial^1\ppartial^{b-1}+ \dots +
\ppartial^{b-1}\ppartial^1+ D\ppartial^{b}\right) = 0.
\]
However, we have already observed that the identity \eqref{twist5}
at stage $b$ is a formal consequence of \eqref{twist4} and, since
the latter has already been established, \eqref{twist5} is now
available at stage $b$, viz.
\[
D\tau^{b+1}=\vartheta_{b+1} - \tau^1 \ppartial^b \colon
\Sigm^{\mathrm c} \longrightarrow \fra g .
\]
Hence
\[
0=D\vartheta_{b+1} - D(\tau^1 \ppartial^b) =D\vartheta_{b+1} +
\tau^1 D\ppartial^b.
\]
Substituting the right-hand side of \eqref{twist4} at stage $b$
for $D\vartheta_{b+1}$, we obtain the identity
\[ 0=\dd^0 \ppartial^{b} + \ppartial^1\ppartial^{b-1}+ \dots +
\ppartial^{b-1}\ppartial^1+ \ppartial^{b}\dd^0 \colon
\Sigm^{\mathrm c} \longrightarrow \Sigm^{\mathrm c},
\]
that is, the identity \eqref{twist22} at stage $b$. This completes
the inductive step. \end{proof}

\section{The proof of the Lie algebra perturbation Lemma}
\label{proofliepert}

Lemma \ref{lem12} entails that the operator $\pppartial$ given by
\eqref{proc2} is a coalgebra perturbation and that the morphism
$\tau$ given by \eqref{proc1} is a Lie algebra twisting cochain.
We will now establish Complement II of the Lie algebra
perturbation Lemma.

The contraction \eqref{2.3} may be obtained in the following way:
Any contraction of chain complexes of the kind \eqref{co}  induces
a filtered contraction
\begin{equation}
   \left(\mathrm T^{\mathrm c}[M]
     \begin{CD}
      \null @>{\mathrm T^{\mathrm c}[\nabla]}>> \null\\[-3.2ex]
      \null @<<{\mathrm T^{\mathrm c}[\pi]}< \null
     \end{CD}
    \mathrm T^{\mathrm c}[N],\mathrm T^{\mathrm c}[h] \right)
\label{2.4}
\end{equation}
of coaugmented differential graded coalgebras. A version thereof
is spelled out as a {\em contraction of bar constructions\/}
already in Theorem 12.1 of \cite{eilmactw}; the filtered
contraction \eqref{2.4} may be found in \cite{gulstatw} (2.2) and
\cite{huebkade} ${\mathrm (2.2.0)}_*$ (the dual filtered
contraction of augmented differential graded algebras being
spelled out explicitly in \cite{huebkade} as ${\mathrm
(2.2.0)}^*$). The differential graded symmetric coalgebras
${\Sigm^{\mathrm c}[M]}$ and ${\Sigm^{\mathrm c}[N]}$ being
differential graded subcoalgebras of ${\mathrm T^{\mathrm c}[M]}$
and ${\mathrm T^{\mathrm c}[N]}$, respectively, the morphisms
${\mathrm T^{\mathrm c}[\nabla]}$ and ${\mathrm T^{\mathrm
c}[\pi]}$ pass to corresponding morphisms ${\Sigm^{\mathrm
c}[\nabla]}$ and ${\Sigm^{\mathrm c}[\pi]}$ respectively, and
${\Sigm^{\mathrm c}[h]}$ arises from ${\mathrm T^{\mathrm c}[h]}$
by {\em symmetrization\/}, so that
\begin{equation}
   \left(\Sigm^{\mathrm c}[M]
     \begin{CD}
      \null @>{\Sigm^{\mathrm c}[\nabla]}>> \null\\[-3.2ex]
      \null @<<{\Sigm^{\mathrm c}[\pi]}< \null
     \end{CD}
    \Sigm^{\mathrm c}[N],\Sigm^{\mathrm c}[h] \right)
\label{2.5}
  \end{equation}
constitutes a filtered contraction of coaugmented differential
graded coalgebras. Alternatively, since $\Sigm^{\mathrm c}$ is a
functor, application of this functor to \eqref{co} yields
\eqref{2.5}. Here $\Sigm^{\mathrm c}[\nabla]$ and $\Sigm^{\mathrm
c}[\pi]$ are morphisms of differential graded coalgebras but,
beware, even though ${\mathrm T^{\mathrm c}[h]}$ is compatible
with the coalgebra structure in the sense that it is a homotopy of
morphisms of differential graded coalgebras, $\Sigm^{\mathrm
c}[h]$ no longer has such a compatibility property in a naive
fashion. Indeed, for differential graded cocommutative coalgebras,
the notion of homotopy is a subtle concept, cf. \cite{schlstas}.
To sum up, application of the functor $\Sigm^{\mathrm c}$ to the
induced contraction \Nsddata {s\fra g}{s\nabla}{s\pi}{sM}{sh}
which arises from \eqref{cont1} by suspension yields the
contraction \eqref{2.3}.

To establish  Complement II of the Lie algebra perturbation Lemma,
we will view the contraction \eqref{2.5} merely as one of filtered
chain complexes, that is, we forget about the coalgebra
structures. As before, we denote by $\partial$ the coalgebra
perturbation on $\Sigm^{\mathrm c}[s\fra g]$ which corresponds to
the graded Lie bracket on $\fra g$, so that the differential on
the {\scs CCE\/}-coalgebra $\mathcal C[\fra g]$ (having
$\Sigm^{\mathrm c}[s\fra g]$ as its underlying coaugmented graded
coalgebra) is given by $d +
\partial$. For intelligibility, we recall the following.

\begin{Lemma}[Ordinary perturbation Lemma]
\label{ordinary} {Let \Nsddata N{\nabla}{\ppi}Mh be a filtered
contraction, let $\partial$ be a perturbation of the differential
on $N$, and let
\begin{align*} \dell &= \sum_{n\geq 0} \ppi\partial (-h\partial)^n\nabla =
\sum_{n\geq 0} \ppi(-\partial h)^n\partial\nabla
\\
\nabla_{\partial}&= \sum_{n\geq 0} (-h\partial)^n\nabla
\\
\ppi_{\partial}&= \sum_{n\geq 0} \ppi(-\partial h)^n
\\
h_{\partial}&=-\sum_{n\geq 0} (-h\partial)^n h =-\sum_{n\geq 0}
h(-\partial h)^n .
\end{align*}
If the filtrations on $M$ and $N$ are complete, these infinite
series converge, $\dell$ is a perturbation of the differential on
$M$ and, when $N_{\partial}$ and $M_{\dell}$ refer to the new
chain complexes,
\begin{equation}
   \left(M_{\dell}
     \begin{CD}
      \null @>{\nabla_{\partial}}>> \null\\[-3.2ex]
      \null @<<{\ppi_{\partial}}< \null
     \end{CD}
   N_{\partial} , h_{\partial} \right)
\label{2.66}
  \end{equation}
constitute a new filtered contraction that is natural in terms of
the given data.}
\end{Lemma}

\begin{proof} See \cite{rbrown} or \cite{gugenhtw}.
\end{proof}

Application of the ordinary perturbation Lemma to the contraction
\eqref{2.5} and the perturbation $\partial$ of the differential
$d$ on $\Sigm^{\mathrm c}[s\fra g]$ yields the perturbation
$\delta$ of the differential $d^0$ on $\Sigm^{\mathrm c}[sM]$ and
the contraction \eqref{2.6a} where the notation $\nnabla$,
$\widetilde \Pi$, $H$ in \eqref{2.6a} corresponds to,
respectively, $\nabla_{\partial}$, $\pi_{\partial}$,
$h_{\partial}$ in \eqref{2.66}. By construction, the composite
\begin{equation*}
\begin{CD}
\Phi\colon\Sigm^{\mathrm c}_{\pppartial}[sM] @>{\overline \tau}>>
\mathcal C[\fra g] @>{\widetilde \Pi}>> \Sigm^{\mathrm
c}_{\delta}[sM]
\end{CD}
\end{equation*}
introduced as \eqref{comp2} above is a morphism of chain complexes
and, modulo the filtrations, as a morphism of the underlying
graded $R$-modules, this composite is the identity. More
precisely, $\Phi$ can be written as an infinite series
\begin{equation}
\Phi = \mathrm{Id} + \Phi^1 + \dots + \Phi^j +\ldots \label{multi}
\end{equation}
such that, for $j \geq 1$, $\Phi^j$ lowers the coaugmentation
filtrations by $j$. Furthermore, the convergence of the series
\eqref{multi} is naive, that is, in each degree, the limit is
achieved after finitely many steps. Consequently $\Phi$ is an
isomorphism of chain complexes. The inverse $\Psi$ of $\Phi$ can
be obtained as the infinite series
\begin{equation}
\Psi = \mathrm{Id} + \Psi^1 + \dots + \Psi^j +\ldots
\label{multi2}
\end{equation}
determined by the requirement $ \Phi\Psi = \mathrm{Id}$ or,
equivalently, $\Psi$ is given by the recursive description
\begin{equation}
\Psi^j +\Phi^1 \Psi^{j-1}   + \dots + \Phi^{j-1}  \Psi^1
+\Phi^j=0,\ j \geq 1,\label{multi3}
\end{equation}
with the convention $\Phi^0=\mathrm{Id}$ and $\Psi^0=\mathrm{Id}$.
Recall from \eqref{Pipert1} and \eqref{Hpert1} that, by
definition,
\begin{align*}
\Pi &= \Psi \widetilde\Pi \colon \mathcal C[\fra g]
\longrightarrow \Sigm^{\mathrm c}_{\pppartial}[sM],
\\
H &= \HH- \HH \overline \tau\,\Pi \colon \mathcal C[\fra g]
\longrightarrow \mathcal C[\fra g].
\end{align*}
By construction,
\[
 \Pi\, \overline \tau
=\mathrm{Id}\colon \Sigm^{\mathrm c}_{\pppartial}[sM]
\longrightarrow \Sigm^{\mathrm c}_{\pppartial}[sM]
\]
and, since \eqref{2.6a} is a contraction of chain complexes,
\begin{align*}
DH &= D(\HH- \HH \overline \tau\,\Pi) =(\mathrm{Id} -\widetilde
\nabla \widetilde \Pi)(\mathrm{Id} -\overline \tau\,\Pi)
\\
&=\mathrm{Id} -\widetilde \nabla \widetilde \Pi  -\overline
\tau\,\Pi + \widetilde \nabla \widetilde \Pi \overline \tau\,\Pi
\\
&= \mathrm{Id} -\widetilde \nabla \Phi \Pi  -\overline \tau\,\Pi +
\widetilde \nabla \widetilde \Pi \overline \tau\,\Pi = \mathrm{Id}
-\overline \tau\,\Pi
\end{align*}
since $\Phi = \widetilde \Pi \overline \tau$ and $\Phi \Pi=
\widetilde \Pi$. Consequently
\[
(d+\partial)H + H (d+\partial) = \mathrm{Id} - \overline \tau
\,\Pi,
\]
that is, $H$ is a chain homotopy between $\mathrm{Id}$ and
$\overline \tau\,\Pi$. Moreover, since \eqref{2.6a} is a
contraction of chain complexes, the side conditions \eqref{side}
hold, that is
\begin{equation}
\Pi  H  = 0, \quad H \,\overline \tau = 0,\quad H H = 0.
\label{sidepert}
\end{equation}
Consequently $\overline {\tau}$, $\Pi$ and $H$ constitute a
contraction of chain complexes of the kind \eqref{cont5}  as
asserted and this contraction is obviously natural in terms of the
data. This establishes Complement II and thus completes the proof
of the Lie algebra perturbation Lemma.


\begin{thebibliography}{19}

\bibitem{rbrown}  R. Brown: The twisted Eilenberg--Zilber
theorem. Celebrazioni Archimedee del Secolo XX, Simposio di
topologia, pp. 33--37 (1964).

\bibitem{cartanse}  H. Cartan: Alg\`ebres
d'Eilenberg--Mac Lane et homotopie. Expos\'es 2--11. S\'eminaire
H. Cartan 1954/55. \'Ecole Normale Sup\'erieure, Paris, 1956.

\bibitem{dold} A. Dold: Zur Homotopietheorie der Kettenkomplexe.
Math. Ann. \textbf{140} (1960), 278--298.

\bibitem{eilmactw} S. Eilenberg and S. Mac Lane: On the
groups ${\mathrm H(\pi,n)}$. I. Ann. of Math. \textbf{58} (1953),
55--106. II. Methods of computation.  Ann. of Math. \textbf{60}
(1954), 49--139.

\bibitem{gugenhtw} V.K.A.M. Gugenheim: On the chain
complex of a fibration.  Illinois J. of Math. \textbf{16} (1972),
398--414.

\bibitem{gugenhth} V.K.A.M. Gugenheim: On a perturbation
theory for the homology of the loop space.  J. of Pure and Applied
Algebra \textbf{25} (1982), 197--205.

\bibitem{gulstatw} V.K.A.M. Gugenheim, L. Lambe, and J.~D.
Stasheff: Perturbation theory in differential homological algebra.
II. Illinois J. of Math. \textbf{35} (1991), 357--373.

\bibitem{homotype} J. Huebschmann: The homotopy type of
$F\Psi^q$. The complex and symplectic cases. In: Applications of
Algebraic $K$-Theory to Algebraic Geometry and Number Theory, Part
II, Proc. of a conf. at Boulder, Colorado, June 12 -- 18, 1983.
Cont. Math. \textbf{55} (1986), 487--518.

\bibitem{perturba} J. Huebschmann: Perturbation theory
and free resolutions for nilpotent groups of class 2.  J. of
Algebra \textbf{126} (1989), 348--399.

\bibitem{cohomolo} J. Huebschmann: Cohomology of
nilpotent groups of class 2.  J. of Algebra \textbf{126} (1989),
400--450.

\bibitem{modpcoho} J. Huebschmann: The mod $p$
cohomology rings of metacyclic groups.  J. of Pure and Applied
Algebra \textbf{60} (1989), 53--105.

\bibitem{intecoho} J. Huebschmann: Cohomology of
metacyclic groups.  Trans. Amer. Math. Soc. \textbf{328} (1991),
1-72.

\bibitem{minimult} J. Huebschmann:
Minimal free multi models for chain algebras. Georgian Math. J.
\textbf{11} (2004), 733--752, {\tt math.AT/0405172}.

\bibitem{koszul} J. Huebschmann: { Homological perturbations,
equivariant cohomology, and Koszul duality}, {\tt
math.AT/0401160}.

\bibitem{koszultw} J. Huebschmann: {Relative homological
algebra, homological perturbations, equivariant de Rham theory,
and Koszul duality}, {\tt math.AT/0401161}.

\bibitem{pertlie2} J. Huebschmann: {The sh-Lie algebra
perturbation Lemma} , {\tt arxiv:0710.2070}.

\bibitem{jimmurra} J. Huebschmann: {Origins and breadth of the theory
of higher homotopies}, {\tt arxiv:0710.2645}.

\bibitem{huebkade} J. Huebschmann and T. Kadeishvili:
Small models for chain algebras.  Math. Z. \textbf{207} (1991),
245--280.

\bibitem{huebstas} J. Huebschmann and J. D. Stasheff:
Formal solution of the master equation via HPT and deformation
theory. Forum Math. \textbf{14} (2002), 847--868.

\bibitem{kontssev}  M.~Kontsevich: {Deformation quantization of
Poisson manifolds},  {\tt math.QA/9709040}.

\bibitem{moorefiv} J. C. Moore: Differential homological
algebra.  Actes, Congr\`es intern. math. Nice (1970),
Gauthiers-Villars, Paris (1971), 335--339.

\bibitem{quilltwo} D. Quillen: Rational homotopy theory.
 Ann. of Math. \textbf{90} (1969), 205--295.

\bibitem{schlstas} M. Schlessinger and J.~D. Stasheff:
Deformation theory and rational homotopy type.  Pub. Math. Sci.
IHES. To appear; new version July 13, 1998.

\bibitem{shih}  W. Shih: Homologie des espaces fibr\'es.
Pub. Math. Sci. IHES \textbf{13} (1962).

\end{thebibliography}
\end{document}